\documentclass[10pt,a4paper]{amsart}
\usepackage{amssymb}
\usepackage[latin1]{inputenc}
\usepackage[english]{babel}
\usepackage{textcomp}
\usepackage{amsthm}
\usepackage{amsmath}
\usepackage{amssymb}
\usepackage[all]{xy}
\usepackage{url}
\newtheorem{teo}{Theorem}[section]
\newtheorem{pro}[teo]{Proposition}
\newtheorem{lem}[teo]{Lemma}
\newtheorem{cor}[teo]{Corollary}
\newtheorem{rem}[teo]{Remark}
\newtheorem{ass}{Assumption}
\newtheorem*{defi}{Definition}

\DeclareFontFamily{OT1}{pzc}{}
\DeclareFontShape{OT1}{pzc}{m}{it}{<-> s * [1.10] pzcmi7t}{}
\DeclareMathAlphabet{\mathpzc}{OT1}{pzc}{m}{it}
\DeclareMathOperator{\Pic}{Pic}
\DeclareMathOperator{\Fix}{Fix}
\DeclareMathOperator{\id}{id}

\DeclareMathOperator{\Hom}{Hom}
\DeclareMathOperator{\divi}{div}
\DeclareMathOperator{\coker}{Coker}
\DeclareMathOperator{\Mo}{Mon}

\newcommand{\pn}{\mathbb{P}}
\newcommand{\rk}{\mathrm{rk}\,}
\newcommand{\sign}{\mathrm{sign}\,}
\newcommand{\hsk}{K3^{\left[2\right]}}
\newcommand{\hskn}{K3^{\left[n\right]}}
\newcommand{\IP}{ \mathbb{P}}
\newcommand{\IC }{ \mathbb{C}}
\newcommand{\IR}{ \mathbb{R}}

\newcommand{\IZ}{\mathbb{Z}}
\newcommand{\IQ}{\mathbb{Q}}
\newcommand{\IN}{\mathbb{N}}

\newcommand{\prf}{{\it Proof.~}}
\newcommand{\ra}{\rightarrow}
\newcommand{\eprf}{\hfill $\square$ \smallskip\par}
\newcommand{\stab}{\mathrm{Stab}}
\begin{document}
 \title[Some remarks on moduli spaces of lattice polarized IHS manifolds]{Some remarks on moduli spaces of lattice polarized holomorphic symplectic manifolds}
\author{Chiara Camere}
\address{Chiara Camere, Universit\`a degli Studi di Milano,
Dipartimento di Matematica ``F. Enriques'',
Via Cesare Saldini 50, 20133 Milano, Italy} 
\email{chiara.camere@unimi.it}
\urladdr{http://www.mat.unimi.it/users/camere/en/index.html}

\begin{abstract}
We construct quasi-projective moduli spaces of $K$-general lattice polarized irreducible holomorphic symplectic manifolds. Moreover, we study their Baily--Borel compactification and investigate a relation between one-dimensional boundary components and equivalence classes of rational Lagrangian fibrations defined on mirror manifolds.
\end{abstract}
\maketitle
\section{Introduction}

This article is intended to develop some aspects of the theory of lattice polarized irreducible holomorphic symplectic manifolds, illustrated in \cite{C4}. In particular, the main focus here will be on how to construct algebraic moduli spaces in this framework, a question that is closely related to the study of the period map, namely of its injectivity.

In the case of $K3$ surfaces, the work of  Dolgachev \cite{dolgachevmirror}, Nikulin \cite{Nikulin} and Pynkham  \cite{Pinkham} shows that it is enough to consider ample polarizations to obtain algebraic moduli spaces, and also gives a precise description of the degenerations on ``bad'' hyperplanes in the period domain. We will see in \S \ref{ample-pol} that the good notion here is less clear geometrically, and this is due to the fact that in higher dimensions birational geometry is highly non-trivial.

The lack of projective models, even in the case of fourfolds of $\hsk$-type, is somewhat of an obstacle to the development of explicit constructions. We will anyways be able to see an example of non-separability phenomenon in \S \ref{ample-pol}.

Once the moduli space under consideration is algebraic, it is natural to look for its compactification, and the smallest one is the Baily--Borel one, which we study in \S \ref{bound-comp}. Classical methods, as developed by Scattone \cite{Scattone}, can be used to construct this compactification; some subtleties arise, with respect to the case of $K3$ surfaces, because of the different arithmetic groups involved. These orthogonal groups are computed explicitly for the two deformation families of manifolds of $\hskn$-type and of generalized Kummer varieties in \S \ref{ortho}.

Finally, we explore in \S \ref{mir-lag} the mirror symmetry introduced in \cite{C4}, and in particular we remark that one-dimensional boundary components of the Baily--Borel compactification give some information on the number, up to birational transformation, of rational Lagrangian fibrations on members of the mirror family.

\subsection*{Acknowledgements}

The author wants to thank Prof. Igor Dolgachev for asking the interesting question about the relation between boundary components and Lagrangian fibrations. This work has greatly benefited from discussions with Prof. S. Boissi\`ere and Prof. A. Sarti. Finally, the author wants to thank Prof. M. Schuett, G. Mongardi and C. Lehn for useful conversations, and A. Garbagnati for all her explanations and for pointing out Nishiyama's paper.
\section{Preliminary notions}

\subsection{Lattices and orthogonal groups}

A lattice $L$ is a free $\IZ$-module endowed with a non-degenerate symmetric bilinear form with integer values; it is even if the associated quadratic form takes values in the even integer numbers. The dual of $L$, usually denoted $L^*$, is $\Hom_{\IZ}(L,\IZ)\cong \lbrace x\in L\otimes\IQ|(x,v)\in\IZ\ \forall v\in L\rbrace$, and the quotient $A_L:=L^*/L$ is the so-called discriminant group. It is a finite group, on which the bilinear form of $L$ induces a finite quadratic form $q_L:A_L\rightarrow \IQ/2\IZ$.
The signature $(r,s)$ of a lattice $L$ is the signature of $L\otimes\IR$; the triple of invariants $(r,s,q_L)$ is called the genus of $L$.

An embedding of a sublattice $j:M\subset L$ is primitive if the quotient $L/M$ is free. A result by Nikulin \cite[Proposition 1.15.1]{Nikulinintegral} gives a powerful technique to classify all possible primitive embeddings of a sublattice up to isometry.

In the sequel, the following lattices will appear: the unique even unimodular hyperbolic lattice of rank two $U$; the even, negative definite lattices $A_k, D_h, E_l$ associated to the Dynkin diagrams of the corresponding type ($k\geq 1$, $h\geq 4$, $l=6,7,8$); and, for $d\equiv -1\ (4)$, the following negative definite lattice 
$$
B_d :=\left(
\begin{array}{cc}
-(d+1)/2&1\\
1&-2\\
\end{array}
\right).
$$
Moreover, $L(t)$ will denote the lattice whose bilinear form is the one on $L$ multiplied by $t\in\IN^\ast$. 

The group $O(L)$ is the group of isometries of $L$. In this paper some of its subgroups will play an important role, namely the
{\it stable orthogonal group}, $\tilde{O}^+(L)$, and its intersection with the special linear group, $ \tilde{SO}^+(L)$, defined respectively as follows:
$$\tilde{O}^+(L)=\left\lbrace g\in O(L)\mid g_{|A_L}=\id, \mathrm{sn}_{\IR}^L(g)=1\right\rbrace,$$
\[\tilde{SO}^+(L):=\left\lbrace g\in \tilde{O}^+(L)\mid \det(g)=1\right\rbrace,\]
where the {\it real spinor norm} $\mathrm{sn}_{\IR}^L:O( L_{\IR})\rightarrow \IR^*/(\IR^*)^2\cong\lbrace \pm 1\rbrace$ is defined as $$\mathrm{sn}_{\IR}^L(g)=\left(-\frac{v_1^2}{2}\right)\cdots\left(-\frac{v_m^2}{2}\right)(\IR^*)^2$$ for $g\in O(L_{\IR})$ factored as a product of reflections $g=\rho_{v_1}\circ\dots\circ\rho_{v_m}$ with $v_i\in L_{\IR}$.
These are well-known arithmetic subgroups of $O(L)$.

\subsection{IHS manifolds}

Irreducible holomorphic symplectic manifolds $X$, often referred to as IHS manifolds for short, are compact K\"ahler manifolds that are simply connected and such that $H^{2,0}(X)$ is spanned by a non-symplectic form $\omega_X$. These are even-dimensional manifolds with $c^1_{\IR}(X)=0$, which is one of the reasons of their importance; the interested reader can refer to \cite{Beauvillec1Nul} and to \cite{Huybrechts}.

Among all smooth complex surfaces, the only example of IHS manifolds is given by $K3$ surfaces. In higher dimensions, there is more than one deformation class, at least two in each dimension $2n$ - the one of Hilbert schemes of $0$-dimensional subschemes of length $n$ of $K3$ surfaces and that of generalized Kummer varieties. There are also two sporadic deformation classes, discovered by O'Grady, in dimension six and ten. As for the present, no more is known.

One peculiar feature of the theory of IHS manifolds is that the group $H^2(X,\IZ)$ carries a quadratic form $q$, the so-called Beauville--Bogomolov--Fujiki quadratic form, which turns it into a non-degenerate lattice of signature $(3,b_2(X)-3)$, conjectured to be even (it is indeed even for all known deformation classes). This lattice is a deformation invariant for a given family; theoretically it does not suffice to determine the deformation class, but in practice it currently identifies one of the known families. If $H^2(X,\IZ)\cong L$ for $L$ an abstract lattice, we talk about IHS manifolds {\it of type $L$}, meaning that we are fixing a deformation class; a {\it marking} of an IHS manifold $X$ of type $L$ is an isometry $\phi:H^2(X,\IZ)\rightarrow L$. 

The main tool to study moduli spaces $\mathcal{M}_L$ of marked IHS manifolds of type $L$ is the period map, $\mathcal{P}:\mathcal{M}_L\rightarrow D_L$, which maps a marked pair $(X,\phi)$ to its period $\left[\phi(H^{2,0}(X))\right]$ inside the {\it period domain}
\[
D_L:=\left\lbrace x\in \IP(L\otimes \IC)|q(x)=0, (x,\overline{x})>0\right\rbrace
\]

The period map is known to be a local isomorphism \cite{Beauvillec1Nul} and surjective \cite{Huybrechts}, also when restricted to one connected component of $\mathcal{M}_L$. A major break-through in the theory has been the proof of the so-called global Torelli theorem, due to Verbitsky \cite{VerbitskyTorelli}, generalizing in a weaker form the one for $K3$ surfaces. Given an IHS manifold $X$ of type $L$, let $\mathcal{C}_X$ and $\mathcal{K}_X$ be respectively the positive and the K\"ahler cone of $X$.

\begin{teo}[Global Torelli Theorem] \cite{VerbitskyTorelli},\cite{HuybrechtsBourbaki},\cite[Theorem 2.2]{MarkmanTorelli} \label{GTT}

Let $\mathcal{M}^+_L$  be a connected component of $\mathcal{M}_L$.
\begin{enumerate}
\item For each $\omega\in D_L$, the fiber $\mathcal{P}^{-1}(\omega)$ consists of pairwise inseparable points.   
\item Let $(X_1,\eta_1)$ and $(X_2,\eta_2)$ be two inseparable points of $\mathcal{M}_L^+$. Then $X_1$ and $X_2$ are bimeromorphic. 
\item The point $(X,\eta)\in\mathcal{M}_L^+$ is Hausdorff  if and only if $\mathcal{C}_X=\mathcal{K}_X$.
 
\end{enumerate}
\end{teo}

In the sequel we will be using also the following Hodge theoretical version of the Torelli theorem.

\begin{teo}\cite[Theorem 1.3]{MarkmanTorelli}\label{HTT}. 
 Let $X$ and $Y$ be two IHS manifolds deformation equivalent one to each other. Then:
 \begin{enumerate}
  \item $X$ and $Y$ are bimeromorphic if and only if there exists a parallel transport operator $f:H^2(X,\IZ)\ra H^2(Y,\IZ)$ that is an isomorphism of integral Hodge structures;
  \item if this is the case, there exists an isomorphism $\tilde{f}:X\ra Y$ inducing $f$ if and only if $f$ preserves a K\"ahler class.
 \end{enumerate}
\end{teo} 

Given $(X,\phi)\in \mathcal{M}_L^+$, let $\Mo^2(L)$ be the group $\phi\circ \Mo^2(X)\circ \phi^{-1}$,  as defined in Markman \cite{MarkmanTorelli}.

\subsection{The K\"ahler cone and the fundamental exceptional chamber}

We denote $\Mo^2_{Hdg}(X)$ the subgroup of monodromies of $X$ that are isomorphisms of integral Hodge structures. Recall that prime exceptional divisors $E$ on $X$ are reduced and irreducible effective divisors such that $E^2<0$; let $\mathcal{P}ex_X$ be the set of all prime exceptional divisors on $X$ and $W_{Exc}(X)$ be the subgroup of  $\Mo^2_{Hdg}(X)$ generated by reflections in prime exceptional divisors.

The {\it fundamental exceptional chamber} $\mathcal{FE}_X$ is the cone $$\lbrace x\in\mathcal{C}_X|(x,E)>0\ \forall E\in\mathcal{P}ex_X\rbrace.$$

\begin{lem}\cite[Proof of Corollary 5.7]{MarkmanTorelli}\label{very-gen-fundam}
Let $X$ and $Y$ be IHS manifolds of type $L$ and take $x\in\mathcal{FE}_{X}$. Then, if there exists a bimeromorphic map $f : X\dashrightarrow Y$, $f_*(x)$ belongs to $\mathcal{FE}_{Y}$.
\end{lem}
\begin{teo}\cite[Theorem 1.6]{MarkmanTorelli}\label{w-exc-decomp}
Let $X$ and $Y$ be IHS manifolds of type $L$ and let $f : H^2(X,\IZ)\rightarrow H^2(Y,\IZ)$ be a parallel transport operator, which is an isomorphism of Hodge structures. Then there exists a
unique element $w \in W_{Exc}(Y)$ and a birational map $g : X\dashrightarrow Y$, such that $f = w \circ g^*$.
The map $g$ is determined uniquely up to composition with an automorphism of $X$, which
acts trivially on $H^2(X,\IZ)$.
\end{teo}

In \cite{MarkmanTorelli}, the author defines the {\it exceptional chambers} of the positive cone $\mathcal{C}_X$ as subsets of the form $g(\mathcal{FE}_X)$, for  $ g\in \Mo^2_{Hdg}(X)$, and the {\it K\"ahler-type chambers} of $\mathcal{C}_X$ as subsets of the form $g(f^*(\mathcal{K}_Y))$,
for  $ g\in \Mo^2_{Hdg}(X)$ and $f : X\dashrightarrow Y$ a bimeromorphic map to an IHS manifold $Y $.

\subsection{Lattice polarized IHS manifolds}

Given $X$ an IHS manifold of type $L$ and $j:M\subset L$ a primitive embedding of a sublattice $M$ of signature $(1,t)$, we defined in \cite{C4} an $M$-polarization of $X$ as a primitive embedding $i:M\rightarrow \Pic(X)$, and a $j$-marking of an $M$-polarized $X$ as a marking $\phi: H^2(X,\IZ)\rightarrow L$ such that $\phi\circ i=j$; in such a case, we say that the pair $(X,\phi)$ is {\it$(M,j)$-polarized}, and it follows from the definition that $\mathcal{P}(X,\phi)\in\pn (M_{\mathbb{C}}^{\perp})$.

The relevant period domain is $$D_{M}=\{\left[\omega \right]\in \pn(N_{\mathbb{C}})\mid
q(\omega)=0,(\omega,\bar{\omega})>0 \},$$ where $N=j(M)^{\perp}$. This has two connected components and each one is a symmetric homogeneous domain of type IV (see \cite{GHSHandbook}).

The construction in \cite{C4} gives a coarse moduli space $\mathcal{M}_{M,j} $ of $(M,j)$-polarized IHS manifolds of type $L$ and a holomorphic map $\mathcal{P}_{M,j}:\mathcal{M}_{M,j}\rightarrow D_M$ that is the restriction of the period map $\mathcal{P}$. Let $D^+_{M}$ be a connected component of $D_M$.

The group $O(L,M)=\{g\in O(L)\mid g(m)=m\ \forall m\in M\}$ acts properly and discontinuously on $D_M$; the group of {\it $(M,j)-$polarized monodromy operators} is $$\Mo^2(M,j):=\{g\in \Mo^2(L)\mid g(m)=m\ \forall m\in M\}=\Mo^2(L)\cap O(L,M) $$ 
It acts on $\mathcal{M}_{M,j}$ via $(X,\phi)\mapsto (X,g\circ\phi)$ for $g\in\Mo^2(M,j)$; let $\Gamma_{M,j}$ be its image in $O(N)$ via the standard restriction map.

\begin{teo}\cite[Proposition 3.4 ]{C4}
Let $\mathcal{M}_{M,j}^+$ be a connected component of $\mathcal{M}_{M,j} $; the period map restricts surjectively to $
\mathcal{P}_{M,j}:\mathcal{M}_{M,j}^+\rightarrow D^+ _{M}
$ and is equivariant with respect to the action of $\Mo^2(M,j)$ and of $\Gamma_{M,j}\subset O(N)$. 
\end{teo}

\section{Ample polarizations }\label{ample-pol}

In most cases the moduli space $\mathcal{M}_{M,j}^+$ is non-separated, and hence non-algebraic; even the first natural attempt to define a notion of ample polarization does not resolve the issue, as we will recall below. In the meanwhile, Amerik and Verbitsky in \cite{AmerikVerbitsky} have developed their theory about MBM classes, and Joumaah, in \cite{Malek}, has employed their results to study the injectivity of the period map in the case of fourfolds of $\hsk$-type with a non-symplectic involution. A generalization of these ideas will lead to the construction of an algebraic moduli space in the more general framework of lattice polarizations.

We consider the positive cone $\left\lbrace x\in M_{\mathbb{R}}\mid (x,x)>0\right\rbrace$ and pick the connected component $C_M$ such that $i(C_M)$ is a subcone of $\mathcal{C}_{X}$, the connected component containing the K\"ahler cone $\mathcal{K}_{X}$ for $(X,\phi)\in\mathcal{M}_L^+$. This is slightly different from what we made in  \cite{C4}, and it is to ensure that the orientation of the positive cone is a deformation invariant (compare with \cite[\S 4]{MarkmanTorelli}).

We defined in \cite{C4} ample lattice polarized IHS manifolds under the following assumption:
\begin{ass}\label{ass1}
Let $(X,\phi)$ be a marked pair of type $L$. There exists a set $\Delta(L)\subset L$ such that the K\"ahler cone $\mathcal{K}_X$ can be described as $$\mathcal{K}_X=\left\lbrace h\in H^{1,1}(X,\mathbb{R})\mid (h,h)>0, (h,\delta)>0\  \forall\delta\in\Delta(X)^+\right\rbrace$$
where $\Delta(X)^+:=\left\lbrace\delta\in \phi^{-1}(\Delta(L))\cap\Pic(X)\mid (\delta,\kappa)>0\right\rbrace$ for $\kappa\in \mathcal{K}_X$ a fixed K\"ahler class.
\end{ass}

Given $\Delta(M):=\Delta(L)\cap M$ and $\delta^{\perp}=\left\lbrace x\in M_{\mathbb{R}}\mid (x,\delta)=0 \right\rbrace$, we fix, for the rest of the paper, a connected component of $C_M\setminus(\cup_{\delta\in\Delta(M)}\delta^{\perp})$ and call it $K$. We say that $(X,\phi)$ is {\it ample $(M,j)$-polarized} if $i(K)$ contains a K\"ahler class.

Recent  work by Amerik and Verbitsky \cite{AmerikVerbitsky} shows that Assumption \ref{ass1} is satisfied. The authors introduce the notion of monodromy birationally minimal class: a non-zero class $x\in H^{1,1}(X,\IQ)$ with $x^2<0$ is said to be \textit{monodromy birationally minimal} (MBM) if there exists $g\in\Mo^2_{Hdg}(X)$ such that $g(x)^{\perp}\subset H^{1,1}(X,\IQ)$ contains a face of a K\"ahler-type chamber.  Let $\Delta(X)$ be the set of all integral MBM classes $x\in H^{1,1}(X,\IZ)$ on $X$; they show the following
\begin{teo}\cite[Theorems 1.19 and 6.2]{AmerikVerbitsky}\label{wallsMBM}
 Let $X$ be an IHS manifold, $\Delta(X)$ as above and $\mathcal{H}:=\cup_{\delta\in\Delta(X)} \delta^{\perp}$.
 Then the K\"ahler cone of $\mathcal{K}_X$ is a connected component of $\mathcal{C}_X\setminus\mathcal{H}$.
\end{teo}
We define $\Delta(L)$ as the set of $\delta\in L$ such that there exists $(X,\phi)\in\mathcal{M}_{M,j}^+$ satisfying $\phi^{-1}(\delta)\in \Delta(X)$.

\begin{lem}\label{assump1}
The set $\Delta(L)$ satisfies Assumption \ref{ass1}, and $\Mo^2(M,j)$ acts on it.
\end{lem}
\prf The result follows immediately by \ref{wallsMBM} and is also proven in \cite[Proposition 7.2]{Malek}.\eprf

From now on we assume that $\rk M\leq \rk L-3$, in order to deal with period domains of dimension at least one. The image of the period map restricted to the moduli space $\mathcal{M}^a_{M,j}$ of ample $(M,j)$-polarized pairs is also computed in Joumaah's PhD thesis, see \cite[Lemma 9.5]{Malek}. 

\begin{lem}\label{lemma1}
Let $(X,\phi)\in \mathcal{M}_{M,j}^+$ be an ample $(M,j)$-polarized IHS manifold. Then $\mathcal{P}_{M,j}(X,\phi)\in D^+ _{M}\setminus\mathcal{D}$, where \[\mathcal{D}:=\bigcup_{\delta\in \Delta(N)}H_{\delta}.\]
\end{lem}
\prf Let $h\in i(K)\cap \mathcal{K}_X$ and suppose that $\mathcal{P}_{M,j}(X,\phi)\in H_{\delta}$ for some $\delta\in \Delta(N)$. This is equivalent to $(\phi(\omega_X),\delta)=0$, hence $\phi^{-1}(\delta)\in\Delta(X)\subset\Pic(X)$. It follows from \cite[Theorem 6.2]{AmerikVerbitsky} that there exists a K\"ahler class $k\in \mathcal{K}_X$ such that $(k,\phi^{-1}(\delta))> 0$, up to replacing $\delta$ with $-\delta$. On the other hand, $(h,\phi^{-1}(\delta))=0$ because $h\in M$ and this gives a contradiction.\eprf

\begin{pro} \label{imageamples}
The restriction of the period map $\mathcal{P}^a_{M,j}:\mathcal{M}^a_{M,j}\rightarrow D^+_M\setminus\mathcal{D} $ is surjective.
\end{pro}
\prf Take $p\in D^+_M\setminus\mathcal{D}$ and let $(X,\phi)$ be an $(M,j)$-polarized pair such that $\mathcal{P}^a_{M,j}(X,\phi)=p$. If $(X,\phi)$ is ample $(M,j)$-polarized, the proof is finished. If not, we construct a birational model of $X$ that belongs to $\mathcal{M}^a_{M,j}$.

Since $p\notin \mathcal{D}$, $i(M)\nsubseteq\phi^{-1}(\delta)^{\perp}$ for any $\delta$ such that $\phi^{-1}(\delta)\in\Delta(X)$, thus it intersects 
$\mathcal{W}:=\mathcal{C}_X\setminus\cup_{\delta\in\Delta(X)}\delta^{\perp}$ in a non-empty open set of maximal dimension equal to $\rk M$. Since $i=\phi^{-1}\circ j$ is an embedding of lattices, it is clear that 
$$i(M\otimes \IR)\cap\mathcal{W}=i(C_M)\cap\mathcal{W}=i(C_M)\setminus\bigcup_{\delta\in\Delta(X)}(\delta^{\perp}\cap i(C_M)).$$
Moreover, $i(\Delta(M))\subset\Delta(X)$, thus $i(K)$ has non-empty intersection with a chamber of $\mathcal{W}$, i.e. a K\"ahler-type chamber $\tilde{\mathcal{K}}$. Hence, there exist $\tilde{X}$ an IHS manifold, a birational morphism $f:X\dashrightarrow\tilde{X}$ and $g$ a Hodge monodromy operator such that $\mathcal{K}_{\tilde{X}}=\sigma^{-1}(\tilde{\mathcal{K}})$ for $\sigma:=g\circ f^*$.

Then the marked pair $(\tilde{X},\tilde{\phi})$ with $\tilde{\phi}:=\phi\circ \sigma$ is ample $(M,j)$-polarized and $\mathcal{P}^a_{M,j}(\tilde{X},\tilde{\phi})=p$ by construction. Indeed, since $p\in D^+_M$ we consider the embedding $\tilde{\imath}:M\subset \Pic(\tilde{X})$ defined by $\tilde{\imath}=(\sigma)^{-1}\circ i$, so 
that $\tilde{\phi}\circ\tilde{\imath}=\phi\circ i=j$. Moreover, $\tilde{\imath}(K)\cap \mathcal{K}_{\tilde{X}}\cong i(K)\cap\tilde{\mathcal{K}} \neq\emptyset$.\eprf

In the sequel, we will generalize to moduli spaces of lattice polarized IHS manifolds the techniques first developed in \cite{Malek} in the case of moduli spaces of pairs of four-folds of $\hsk$-type with non-symplectic involutions, and then generalized to complex ball quotients in a paper of the author joint with S. Boissi\`ere and A. Sarti \cite{BCSBallQuots}.

We need to recall the following facts from \cite{Malek}: the author defines the set $\Delta'_M\subset L$ as the set of elements $\delta\in L$ such that $\sign(M\cap\delta^{\perp})=(1,\rk(M)-2)$ and $\sign(N\cap \delta^{\perp})=(2, \rk(N)-3)$, then he gives the following characterization.

\begin{lem}{\cite[Lemma 7.6]{Malek}}\label{Malek}
The following are equivalent:
\begin{enumerate}
\item $\delta\in\Delta'_M$;
\item $\delta\notin M\cup N$, $D_M\cap H_{\delta}\neq\emptyset$ and $j(C_M)\cap\delta^{\perp}\neq\emptyset$;
\item given the decomposition $\delta=\delta_M+\delta_N$ with $\delta_M\in M_{\IQ}$ and $\delta_N\in N_{\IQ}$, we have $\delta_M^2<0$ and $\delta_N^2<0$.
\end{enumerate}
\end{lem}

We define $\Delta'(K)$ as the set of $\delta\in\Delta'_M$ such that $\delta^{\perp}\cap j(K)\neq\emptyset$.

\begin{pro}\label{fiber-K}
 Given $p\in D^+_M\setminus\mathcal{D}$, there is a bijection between $(\mathcal{P}^a_{M,j})^{-1}(p)$ and the set of connected components in the chamber decomposition 
 \begin{equation}\label{K-decomp}
  K\setminus (\bigcup_{\delta\in\Delta'(K)\cap p^{\perp}}\delta^{\perp}).
 \end{equation}

\end{pro}
\prf
Given $(X,\phi)\in(\mathcal{P}^a_{M,j})^{-1}(p)$,  we define $\Xi(X,\phi):=\phi(\mathcal{K}_X)\cap K$. It is indeed a chamber of (\ref{K-decomp}): $\mathcal{K}_X$ is a connected component of $\mathcal{C}_X\setminus (\bigcup_{\delta\in\Delta(X)}\delta^{\perp})$ and $\left(\mathcal{C}_X\setminus(\bigcup_{\delta\in\Delta(X)}\delta^{\perp})\right) \cap i(K)=i(K)\setminus (\bigcup_{\delta\in\Delta(X)}\delta^{\perp})$; on the other hand, if $\delta^{\perp}\cap i(K)\neq\emptyset$ for $\delta\in \Delta(X)$, automatically $\delta$ satisfies (ii) of Lemma \ref{Malek}, hence $\delta\in \phi^{-1}(\Delta'(K))\cap \Delta(X)$.

Let us show the injectivity of $\beta$. Given $(X,\phi_X),(Y,\phi_Y)\in (\mathcal{P}^a_{M,j})^{-1}(p)$ two birational models, if $\phi_X(\mathcal{K}_X)\cap K\neq \phi_Y(\mathcal{K}_Y)\cap K$, $\mathcal{K}_X$ and $\phi_X(\phi_Y^{-1}(\mathcal{K}_Y))$ are two distinct K\"ahler-type chambers in $\mathcal{C}_X$, and 
hence $(X,\phi_X)$ and $(Y,\phi_Y)$ cannot be isomorphic by the global Torelli theorem \ref{GTT}.

Finally, fix $\overline{K}$ a chamber of (\ref{K-decomp}) and let $(X,\phi)\in(\mathcal{P}^a_{M,j})^{-1}(p)$ be an element of the fibre. Let $\overline{\Delta}$ be the subset of classes $\delta$ in $\Delta(X)$ such that $(\delta,k)>0$ for all $k\in i(\overline{K})$. Let $\overline{\mathcal{K}}$ be one of the connected components of $\mathcal{C}_X\setminus (\bigcup_{\delta\in\Delta(X)}\delta^{\perp})$ such that $(h,\delta)>0$ for all $h\in \overline{\mathcal{K}}$ and for all $\delta\in \overline{\Delta}$. It is by definition a K\"ahler-type chamber, hence there exist $g\in\Mo^2_{Hdg}(X)$ and a birational morphism $f:X\dashrightarrow \overline{X}$ such that $\overline{\mathcal{K}}=g(f^*\mathcal{K}_{\overline{X}})$. Then, by the proof of Proposition \ref{imageamples}, $(\overline{X}, \phi\circ g\circ f^*)\in(\mathcal{P}^a_{M,j})^{-1}(p)$, and $\Xi(\overline{X}, \phi\circ g\circ f^*)=\overline{K}$.
\eprf

Recall that $(X,\phi)$ is said to be {\it strictly $(M,j)$-polarized} if the embedding $i$ gives an isomorphism $M\cong \Pic(X)$ of hyperbolic lattices.
It is clear that, when elements in $(\mathcal{P}^a_{M,j})^{-1}(p)$ are strictly $(M,j)$-polarized, $\Delta'(K)\cap p^{\perp}=\emptyset$ and (\ref{K-decomp}) consists in one connected component (see \cite[Lemma 3.3]{C4}), so that by Proposition \ref{fiber-K} we recover the fact that the restriction of the period map $\mathcal{P}^a_{M,j}$ to the moduli space of strictly $(M,j)$-polarized manifolds of type $L$ is injective. On the other hand, whilst requiring $\Pic(X)\cong M$ is too restrictive to get an algebraic moduli space (see \cite[Theorem 3.9]{C4}), we will see that the natural condition to impose is the existence of only one connected component in (\ref{K-decomp}).

In particular, we introduce the following notion (see also \cite{BCSBallQuots}).

\begin{defi}
Given the choice of a chamber $K$ of $C_M\setminus(\cup_{\delta\in\Delta(M)}\delta^{\perp})$, an $(M,j)$-polarized pair $(X,\phi)$ is $K$-general if $i(K)=\mathcal{K}_X\cap i(C_M)$.
\end{defi}

Obviously, an $(M,j)$-polarization which is $K$-general is also ample.

\begin{lem}\label{image1}
Let $(X,\phi)\in \mathcal{M}_{M,j}^+$ be a $K$-general $(M,j)$-polarized IHS manifold. Then $\mathcal{P}_{M,j}(X,\phi)\in D^+ _{M}\setminus(\mathcal{D}\cup\mathcal{D}'_K)$, with $\mathcal{D}$ as in Lemma  \ref{lemma1} and \[\mathcal{D}'_K:=\bigcup_{\delta\in \Delta'(K)}H_{\delta}.\]
\end{lem}

\prf We argue by contradiction. If $\mathcal{P}_{M,j}(X,\phi)\in H_{\delta}$ for $\delta\in\Delta(N)$, Lemma \ref{lemma1} implies that $(X,\phi)$ is not even ample $(M,j)$-polarized. If $\mathcal{P}_{M,j}(X,\phi)\in H_{\delta}$ for $\delta\in\Delta'(K)$, we get that $\phi^{-1}(\delta)\in\Delta(X)\subset\Pic(X)$; moreover, there exists $h\in\delta^{\perp}\cap j(K)$. Hence $\phi^{-1}(h)\in\phi^{-1}(\delta) ^{\perp}\cap i(K) = \phi^{-1}(\delta) ^{\perp}\cap\mathcal{K}_X\cap i(C_M)$; in particular, $\phi^{-1}(h)\in\phi^{-1}(\delta) ^{\perp}\cap\mathcal{K}_X$ in contradiction with Theorem \ref{wallsMBM}.\eprf

\begin{rem}
One can show that the isometry class of $i(K)$ is a complete deformation invariant for ample $(M,j)$-polarized pairs $(X,\phi)$. The proof is easily adapted from \cite[\S 9]{Malek}.
\end{rem}

Let $\mathcal{M}_{K,j}^+$ be a connected component of the moduli space of $K$-general $(M,j)$-polarized manifolds inside $\mathcal{M}_{M,j}^+$.

\begin{pro}\label{pksurjective}
The restriction of the period map $\mathcal{P}_{K,j}^+:\mathcal{M}_{K,j}^+\rightarrow D^0_M:=D^+ _{M}\setminus(\mathcal{D}\cup\mathcal{D}'_K)$ is surjective and $\Mo^2(M,j)$-equivariant.
\end{pro}
\prf 
In Lemma \ref{image1} we have shown that the image $\mathcal{P}_{K,j}^+(\mathcal{M}_{K,j}^+)$ is contained in $D^0 _{M}$. Vice versa, if $\omega\in D^0 _{M}$, let $(X,\phi)\in(\mathcal{P}^a_{M,j})^{-1}(\omega)$; by Proposition \ref{imageamples}, 
it is an ample $(M,j)$-polarized pair. In fact, $(X,\phi)$ is also $K$-general: if there was a class $k\in i(K)\setminus\left(\mathcal{K}_X\cap i(C_M)\right)$, there would exist $\delta\in\Delta(X)$ such that $(k,\delta)\leq 0$ and $\delta^{\perp}\cap i(K)\neq 0$, and this would imply $\omega\in H_{\phi(\delta)}$.

The equivariance is obvious.\eprf

\begin{teo}
The period map induces a bijection
\[\mathcal{P}_{K,j}^+:\mathcal{M}_{K,j}^+/\Mo^2(M,j)\rightarrow D^0 _{M}/\Gamma_{M,j}\]
\end{teo}
\prf We only need to show that the period map is injective, because of Proposition \ref{pksurjective}. 
If $(X,\phi),(X',\phi')\in (\mathcal{P}_{K,j}^+)^{-1}(\pi)$ for $\pi\in D^0_M$, 
then $\phi^{-1}\circ\phi ':H^2(X',\IZ)\rightarrow H^2(X,\IZ)$ is a Hodge parallel transport operator. Moreover, $\phi^{-1}(\phi '(\mathcal{K}_{X'}))\cap i(C_M)=\phi^{-1}(j(K))=i(K)\supset\mathcal{K}_X\cap i(C_M)$. Since $\phi^{-1}\circ\phi '$ sends a K\"ahler class on $X'$ to a K\"ahler class on $X$,  by Theorem \ref{GTT} there exists an isomorphism $f:X\dashrightarrow X'$ such that $f^*=\phi^{-1}\circ\phi '$.
\eprf

\begin{cor}
If $\Gamma_{M,j}\subset O(N)$ is an arithmetic subgroup, the moduli space $\mathcal{M}_{K,j}^+/\Mo^2(M,j)$ is a quasi-projective variety of dimension $h^{1,1}-\rk M-2$.
\end{cor}
\prf The quotient  $(\mathcal{D}\cup\mathcal{D}'_K)/\Gamma_{M,j}$ is Zariski-closed: the proof of \cite[Lemma 6.4]{BCSBallQuots} shows that $\mathcal{D}\cup\mathcal{D}'_K$ is locally finite with respect to the action of any subgroup of finite index in $O(N)$ that acts properly discontinuously on $D_M$ and leaves $\Delta(N)$ and $\Delta'(K)$ invariant.\eprf 

\begin{rem}
Also in the case of $K3$ surfaces, the notion of $K$-general $K3$ surfaces is stronger than that of ample $M$-polarized $K3$ surface, so that our result is weaker than Corollary 3.2 in \cite{dolgachevmirror}.
\end{rem}

{\bf Example.} In \cite[\S 7]{Hassett-Tschinkel-flops} the authors compute the ample cone of the Fano variety of a cubic fourfold containing a smooth cubic scroll and of its  Mukai flops. The cubics they start from belong to the Hassett divisor $\mathcal{C}_{12}$ and their Fano varieties are naturally $(M,j)$-polarized, with $M=\IZ g\oplus\IZ\tau\cong\langle 6\rangle\oplus\langle -4\rangle$ (the reader should note that we use a different basis for $M$) and $j$ the unique embedding such that $N=U\oplus E_8^{\oplus 2}\oplus B_3\oplus \langle 4\rangle$. The authors compute that the positive cone $C_M$ is spanned by $(-2+\sqrt{6})g+(3-\sqrt{6})\tau$ and $(2+\sqrt{6})g-(3+\sqrt{6})\tau$.

More precisely, given $X$ a general smooth cubic fourfold containing a smooth cubic scroll $T$, its Fano variety $F$ will be strictly $(M,j)$-polarized and it will contain two Lagrangian planes $P$ and $P^{\vee}$; by flopping each of these planes they obtain two non-isomorphic birational models $F_1$ and $F_1^{\vee}$. The authors then compute the nef cones $K_1$, $K_2$ and $K_3$ of the different models $F$, $F_1$ and $F_1^{\vee}$, which are respectively generated by $\alpha_1=4g+3\tau$ and $\alpha_1^{\vee}=4g-3\tau$, by $\alpha_1$ and $\alpha_2=8g+9\tau$, and finally by $\alpha_1^{\vee}$ and $\alpha_2^{\vee}=8g-9\tau$.

Moreover, in \cite[Theorem 31]{Hassett-Tschinkel-flops} they show that these are the only three birational models of strictly $(M,j)$-polarized fourfolds of $\hsk$ -type up to isomorphism, and they are examples of $K_i$-general $(M,j)$-polarized Fano varieties for $i=1,2,3$.

Now choose the connected component of $\mathcal{M}_{K_1,j}$ containing $(F,\phi)$ with $\phi$ an appropriate marking, and call it  $\mathcal{M}_{K_1,j}^+$. Let $a,b$ and $z$ be respectively the generators of the two summands $B_3$ and $\langle 4\rangle$ in $N$; take $\delta_M:=-\frac{1}{3}g+\frac{3}{2}\tau$, $\delta_N:=\frac{14}{3}a+\frac{8}{3}b+\frac{9}{2}z$ and $\delta:=\delta_M+\delta_N$. We can easily show that $\delta\in L$: let $v_i,w_i$ be a basis of a $U$ summand inside $L$ for $i=1,2$ and $e$ a generator of $\langle -2\rangle$; then the embedding $j$ can be given by $g=2v_1+14w_1-5e$, $\tau=v_2-2w_2$, while $j^{\perp}$ by $a=v_1+3w_1-2e$, $b=e-5w_1$ and $z=v_2+2w_2$. In these coordinates, we get $\delta=4(v_1-w_1)+6(v_2+w_2)-5e\in L$, hence we obtain a class of square $-10$ and divisibility two in $L$ but not in $M\cup N$; moreover, $\delta^2_M<0$ and $\delta^2_N<0$.

On the other hand, the intersection $\delta^{\perp}\cap K_1$ is non-empty, e.g. it contains $\beta_1:=-3g+\tau$, so that $\delta\in \Delta(K_1)$. If we consider a general $\omega\in D^+_M\cap H_\delta$, then in $(\mathcal{P}_{K_1,j}^+)^{-1}(\omega)$ there are two inseparable points $(Y_1,\eta_1)$ and $(Y_2,\eta_2)$ such that $\eta_i(\mathcal{K}_{Y_i})\cap j(C_M)$, $i=1,2$, are respectively the two subchambers of $j(K_1)\setminus\delta^{\perp}$. These two points are clearly not $K_1$-general.

\section{Orthogonal groups}\label{ortho}

In this section we focus on the group $\Gamma_{M,j}$: we have shown in \cite[Proposition 3.5]{C4} that it is arithmetic whenever $\Mo^2(L)$ contains the subgroup $\tilde{O}^+(L)$; however, this criterion can be used only in the case of fourfolds of $\hsk$-type. Here we give a more general result and use it to deal with varieties of $\hskn$-type and of generalized Kummer type.

We denote
$\tilde{SO}^+(L):=\left\lbrace g\in \tilde{O}^+(L)\mid \det(g)=1\right\rbrace.$

\begin{pro}\label{finiteindex}
 If $\Mo^2(L)\supset \tilde{SO}^+(L)$ then $\Gamma_{M,j}$ is an arithmetic subgroup of $O(N)$.
\end{pro}
\prf We show that $H:=\tilde{SO}^+(N)$ is contained in $\Gamma_{M,j}$: since $H$ is of finite index in $O(N)$, this ends the proof. Take $g\in H$; the isometry $f$ of $L$ induced by $\id_M\oplus g$ is an element of $\tilde{O}^+(L)$ by the proof of \cite[Proposition 3.5]{C4}. In particular, one can decompose the extension of $g$ by linearity to $N_{\IR}$ as $g_{\IR}=\rho_{v_1}\circ\dots\circ\rho_{v_m}$ with $v_1,\dots,v_m\in N_{\IR}$, and hence also $f_{\IR}=\rho_{v_1}\circ\dots\circ\rho_{v_m}$ in $O(L_{\IR})$. This implies that $\det(g_{\IR})=\det(f_{\IR})=\prod\det(\rho_{v_i})$ and $\det(g)=\det(f)=1$.\eprf

It follows from Proposition \ref{finiteindex} and from results by Markman \cite{MarkmanTorelli},\cite{MarkmanMehrotra} and Mongardi \cite{MongardiMonodromy} that $\Gamma_{M,j}$ is an arithmetic subgroup of $O(N)$ also for $\hskn$-type with $n\geq 3$ and for generalized Kummer manifolds. Indeed, $\Mo^2(L)$ is respectively isomorphic to 
 $$\hat{O}^+(L):=\left\lbrace g\in O^+(L)\mid g_{|A_L}=\pm \id_{A_L}\right\rbrace$$ if $L=U^{\oplus 3}\oplus E_8^{\oplus 2}\oplus \langle -2(n-1)\rangle$, $n\geq 3$, and to
 $$\hat{SO}^+(L):=\left\lbrace g\in O^+(L)\mid (g_{|A_L}= \id_{A_L},\ \det(g)=1)\ \mathrm{or}\ (g_{|A_L}= -\id_{A_L},\ \det(g)=-1)\right\rbrace$$ if $L=U^{\oplus 3}\oplus \langle -2(n+1)\rangle$.
\begin{cor}
The moduli space $\mathcal{M}_{K,j}^+/\Mo^2(M,j)$ of $K$-general $(M,j)$-polarized pairs of $\hskn$-type or of generalized Kummer type is a quasi-projective variety of dimension $h^{1,1}-\rk M-2$ for any positive integer $n$.
\end{cor}

\begin{teo}
 For $L=U^{\oplus 3}\oplus E_8^{\oplus 2}\oplus \langle -2(n-1)\rangle$ with $n\geq 3$ we have that $$\Gamma_{M,j}= \left\lbrace\begin{array}{l}
                                                                                                            \hat{O}^+(N)\ \mathrm{if}\ A_L\subset A_N\\
                                                                                                            \tilde{O}^+(N)\ \mathrm{otherwise}
                                                                                                           \end{array}\right..
$$
 
 If $L=U^{\oplus 3}\oplus \langle -2(n+1)\rangle$, $\Gamma_{M,j}=\left\lbrace\begin{array}{l}\hat{SO}^+(N)\ \mathrm{if}\ A_L\subset A_N\\
                                                                    \tilde{SO}^+(N)\ \mathrm{otherwise}
                                                                   \end{array}\right.
$.
\end{teo}
\prf This follows from Proposition \ref{finiteindex}, \cite[Proposition 3.5]{C4} and the following fact: given $f\in \Gamma_{M,j}$, then $f$ is the isometry induced on $L$ by $\id_M\oplus f_{|N}$ and we know that $f_{|A_M}=\id_{A_M}$. Hence $f_{|A_L}=-\id_{A_L}$ can happen only if $A_L\subset A_N \subset A_M\oplus A_N$. On the other hand, the same proof as before shows that $\det f=\det f_{|N}$.\eprf

\section{Compactifications}\label{compact}

In this section we study the smallest possible compactification of $D_M^0/\Gamma_{M,j}$, which is the Baily--Borel compactification. Finding a geometrical description of the elements in the boundary is out of the scope of this paper; here, we limit ourselves to specialize the theory in our setting, to give some examples and to discuss an interesting relation, via mirror symmetry, with Lagrangian fibrations.

\subsection{The Baily--Borel compactification}

We introduce now the Baily--Borel compactification of the period domain, defined as its closure in the Harish-Chandra embedding (see \cite{GHSHandbook} for a nice survey of the topic). Let $D^*_M$ be the closure of $D^+_M$ inside the quadric $\{\left[\omega \right]\in \pn(N_{\mathbb{C}})\mid
(\omega,\omega)=0\}.$ 

A {\it boundary component} is of the form $\IP(I_{\IC})\cap D_{M}^*$ for some isotropic subspace $I\subset N_{\IR}$ of dimension one or two; it is called \textit{rational} if the corresponding $I$ can be defined over $\IQ$. In particular, $0$-dimensional rational boundary components of $D^+_{M}$ are in bijection with primitive isotropic elements of $N$.

Let  $\mathcal{RB}$ be the set of rational boundary components of $D^+_{M}$, and, for $F\in\mathcal{RB}$, let $N(F)=\{g\in \Gamma_{M,j}\mid g(F)=F\}$ be  its stabilizer. The Baily--Borel compactification is defined as
\begin{displaymath}
 \overline{D^+_{M}/\Gamma_{M,j}}=D^+_{M}/\Gamma_{M,j}\coprod \left(\coprod_{F\in\mathcal{RB}/\Gamma_{M,j}}F/N(F)\right)
\end{displaymath}
which is a normal projective algebraic variety.

Hence we need to compute  $\Gamma_{M,j}$-orbits of isotropic sublattices of $N$ of rank one or two, in order to obtain the exact number of boundary components in $\overline{D^+_{M}/\Gamma_{M,j}}$.

\subsection{Boundary components and mirror Lagrangian fibrations}\label{mir-lag}

In \cite{dolgachevmirror} the author highlights an interesting relation between one-dimensional boundary components of the Baily--Borel compactification and elliptic fibrations on mirror $K3$ surfaces.

Indeed, fix a primitive embedding $j:M\subset L$ and consider the moduli space $\mathcal{M}^+_{M,j}$. Suppose that there exists a primitive embedding $U\hookrightarrow N$, and let $f,f'$ be a basis of $U$ of isotropic vectors; assume furthermore that $N\cong U\oplus \check{M}$, with $\check{M}\cong (\IZ f)^{\perp}/\IZ f$. Assume for the rest of this section that $U$ admits a unique embedding inside $N$ up to the action of $O(N)$.

The choice of an $O(N)$-orbit of a primitive isotropic $f\in N$ thus implies the choice of a splitting $N\cong U\oplus \check{M}$ and hence of a mirror moduli space $\mathcal{M}_{\check{M},\check{\jmath}}^+$, as illustrated in \cite{dolgachevmirror} and in \cite[\S 4.2]{C4}; on the other hand, the choice of a $\Gamma_{M,j}$-orbit of $f\in N$also identifies a boundary point $F$ of $\overline{D^+_{M}/\Gamma_{M,j}}$, as we have briefly discussed above.

Fix now a  primitive isotropic $f\in N$; we look for isotropic sublattices $S\subset N$ of rank two containing $f$, i.e. for boundary curves passing through the point $F$ identified by $f$. Let $g\in S$ be an element such that $f,g$ is a basis of $S$: $g$ is isotropic and we have $(f,g)=0$. This implies that $g$ is an isotropic element in $\check{M}$; in other words, orbits of isotropic elements of $\check{M}$ with respect to the action of the stabilizer $\stab_{\Gamma_{M,j}}(f)$ are in bijection with rational boundary curves passing through $F$. Let $Z_g$ be the boundary curve passing through $F$ and identified by $g\in\check{M}$.

Given a very general $(Y,\eta)\in\mathcal{M}^+_{\check{M},\check{\jmath}}$ strictly polarized and a primitive isotropic element $g\in\check{M}$, the isomorphism $\check{\imath}:\check{M}\cong\Pic(Y)$ gives a primitive isotropic element $\check{\imath}(g)=c_1(\mathcal{L})$ for $\mathcal{L}$ a line bundle on $Y$ . Conjecturally, a primitive algebraic class $c_1(\mathcal{L})$ in the boundary of the closure of the positive cone $\mathcal{C}_Y$ induces a {\it rational Lagrangian fibration} on $Y$, i.e. there exist a birational map $\psi:Y\dashrightarrow\tilde{Y}$, a Hodge monodromy operator $w\in\Mo^2(Y)$ and a line bundle $\tilde{\mathcal{L}}$ on $\tilde{Y}$ such that $\tilde{\mathcal{L}}$ induces a Lagrangian fibration $\epsilon: \tilde{Y}\rightarrow \IP^n$, with $w(c_1(\mathcal{L}))=\psi^*(\epsilon^*(c_1(\mathcal{O}_{\IP^n}(1))))$.

This conjecture is motivated from what holds true in the case of elliptic fibrations on $K3$ surfaces. Nowadays it is known for two of the four known deformation classes of IHS manifolds, varieties of $\hskn$-type and generalized Kummer varieties, proven respectively by Markman \cite[Theorem 6.3]{MarkmanLagrangian} and Bayer and Macr\`i \cite{BayerMacri} in the Hilbert scheme case, and by Yoshioka, \cite{Yoshioka}, for generalized Kummers varieties; finally it was extended to all their deformations by Matsushita.

\begin{teo}\cite[Corollary 1.1]{MatsushitaIsotropic}. \label{LagrRatFibr}
Let $X$ be an IHS manifold of $\hskn$-type  or of generalized Kummer type of dimension $2n$. We also let $\mathcal{L}$ be a line bundle on $X$ which is not trivial, isotropic with respect to the Beauville--Bogomolov--Fujiki quadratic form on $H^2(X,\IC)$ and such that $c_1(\mathcal{L})$ belongs to $\overline{\mathcal{BK}_X}$. Then $\mathcal{L}$ defines a rational Lagrangian fibration over $\IP^n$.
\end{teo}

\begin{teo}\label{mirrorLag}
Let $\mathcal{M}^+_{M,j}$ be the moduli space of $(M,j)$-polarized varieties of $\hskn$-type or of $(M,j)$-polarized varieties of Kummer type of dimension $2n$. Let $\overline{D^+_{M}/\Gamma_{M,j}}$ be the Baily--Borel compactification of the corresponding period domain. Let $F$ be a boundary point, $N=U\oplus \check{M}$ be the decomposition induced by $F$, as above, and let $\mathcal{M}_{\check{M},\check{\jmath}}^+$ be a connected component of the mirror moduli space. 

Every boundary curve $Z_g$, for $g\in\check{M}$ primitive isotropic, passing through $F$ corresponds to a rational Lagrangian fibration on every strictly $(\check{M},\check{\jmath})$-polarized variety $(Y,\eta)\in\mathcal{M}_{\check{M},\check{\jmath}}^+$.
\end{teo}
\prf
As above, let  $S$ be  a rank two isotropic sublattice of $N$  corresponding to $Z_g$; we can choose as a basis of $S$ the one given by $f\in U$, $g\in\check{M}$. Given  $(Y,\eta)\in\mathcal{M}_{\check{M},\check{\jmath}}^+$ a strictly $(\check{M},\check{\jmath})$-polarized variety, by definition we have an isomorphism $\check{\imath}: \check{M}\cong \Pic(Y)$ such that $\eta\circ\check{\imath}=\check{\jmath}$; the K\"ahler-type chambers of $Y$ are the connected components of $\check{\imath}(C_{\check{M}}\setminus \cup _{\delta\in\Delta(\check{M})}\delta^{\perp})$ (see \cite[Lemma 3.3]{C4}). Hence, there exist a Hodge monodromy operator $w\in\Mo^2(Y)$ and a line bundle $\mathcal{L}$ on $Y$ such that the element $w(\check{\imath}(g))$ is, up to $\pm 1$, an element $c_1(\mathcal{L})$ in $\Pic(Y)$, isotropic with respect to the Beauville--Bogomolov--Fujiki quadratic form and in the boundary of $\mathcal{BK}_Y$. It follows from Theorem \ref{LagrRatFibr} that $\mathcal{L}$ induces a rational Lagrangian fibration on $Y$.

Consider now two isotropic sublattices $S_1,S_2$ of rank two in $N$, and denote $f_k,g_k$ the basis of $S_k$, for $k=1,2$. Suppose that $S_1$ and $S_2$ are in the same $\Gamma_{M,j}$-orbit, and more precisely that there exists $\gamma \in \Gamma_{M,j}$ such that $\gamma(f_1)=f_2$, $\gamma(g_1)=g_2$, so that $F_1=F_2$ and $Z_{g_1}=Z_{g_2}$ identify the same rational boundary components in $\mathcal{RB}/\Gamma_{M,j}$. Each of the $f_k$ induces a decomposition $N\cong U\oplus \check{M}_k$, for $k=1,2$, such that $g_k\in\check{M}_k$; $\gamma$ induces an isomorphism $\check{M}_1\cong\check{M}_2$. The sublattice $\check{M}_k$ is embedded in $\mathcal{L}$ via $\check{\jmath}_k$ defined as the composition $\check{M}_k\hookrightarrow N\hookrightarrow L$; thus we get $\check{\jmath}_1=\check{\jmath}_2\circ \gamma_{|\check{M}_1}$.

Given $(Y,\eta)$ strictly $(\check{M}_1,\check{\jmath}_1)$-polarized, we remark that $(Y,\eta)$ is also strictly $(\check{M}_2,\check{\jmath}_2)$-polarized: indeed, if $\check{\imath}_1:= \eta^{-1}\circ \check{\jmath}_1:\check{M}_1\hookrightarrow\Pic(Y)$, then it is enough to define $\check{\imath}_2:=\check{\imath}_1 \circ\gamma_{|\check{M}_2}^{-1} $. We have $\check{\imath}_1(g_1)=\check{\imath}_2(\gamma(g_1))=\check{\imath}_2(g_2)$, hence they identify the same algebraic element $\mathcal{L}$ in $\Pic(Y)$. 
\eprf

We are interested in studying equivalence classes of birational rational Lagrangian fibrations.
 
\begin{defi}
Given two line bundles $\mathcal{L}_1$ and $\mathcal{L}_2$ on a projective IHS manifold $Y$, let, for $i=1,2$, $\psi_i:Y\dashrightarrow\tilde{Y_i}$ be the birational maps,  $w_i\in\Mo^2(Y)$ be the Hodge monodromy operators and $\tilde{\mathcal{L}}_i$ be the line bundles on $\tilde{Y}_i$, such that $\tilde{\mathcal{L}}_i$ induces a Lagrangian fibration $\epsilon_i: \tilde{Y}_i\rightarrow \IP^n$, with $w_i(c_1(\mathcal{L}_i))=\psi_i^*(c_1(\tilde{\mathcal{L}}_i))$.

We say that $\mathcal{L}_1$ and $\mathcal{L}_2$ induce {\it birational rational Lagrangian fibrations} on $Y$ if there exists a birational map $\sigma:\tilde{Y}_1\dashrightarrow \tilde{Y}_2$ such that $\sigma^*\tilde{\mathcal{L}}_2=\tilde{\mathcal{L}}_1$.
\end{defi}

In particular,  $\epsilon_2\circ \sigma=\epsilon_1$.

\begin{pro}\label{Hdg-mon-Y}

With the same notation as in the Definition, $\mathcal{L}_1$ and $\mathcal{L}_2$ induce birational rational Lagrangian fibrations on $Y$ if and only if there exists a Hodge monodromy $\gamma\in\Mo^2_{Hdg}(Y)$ such that $\gamma(c_1(\mathcal{L}_2))=c_1(\mathcal{L}_1)$.
 
\end{pro}
\prf
We use the same notations as in the definition above. If $\mathcal{L}_1$ and $\mathcal{L}_2$ induce birational rational Lagrangian fibrations on $Y$, the isometry $\gamma:= w_1^{-1}\circ\psi_1^*\circ\sigma^*\circ (\psi_2)_*\circ w_2$ satisfies $\gamma(c_1(\mathcal{L}_2))=c_1(\mathcal{L}_1)$ and it is a Hodge monodromy on $Y$ because it is a composition of parallel transport operators which are isomorphisms of Hodge structures.

Conversely, we define $\overline{f}:=(\psi_2)_*\circ w_2\circ\gamma\circ w_1^{-1}\circ \psi_1^*$; it is a parallel transport operator and an isomorphism of Hodge structures. By Theorem \ref{w-exc-decomp} applied to $\overline{f}^{-1}$, there exist $w\in W_{Exc}(\tilde{Y}_1)$ and a birational morphism $h:\tilde{Y}_1\dashrightarrow\tilde{Y}_2$ such that $\overline{f}=h_*\circ w$. Let $\alpha_i:=c_1(\tilde{\mathcal{L}}_i)$ for $i=1,2$; we have $w(\alpha_1)=h^*(\alpha_2)$. On the other hand, both $\alpha_1$ and $\alpha_2$ are nef, since the corresponding line bundles induce a Lagrangian fibration, hence they belong to the closures of the fundamental exceptional chambers $\overline{\mathcal{FE}_{\tilde{Y}_i}}$. 
 
It follows from Proposition \ref{very-gen-fundam} that $h^*(\overline{\mathcal{FE}_{\tilde{Y}_2}})=\overline{\mathcal{FE}_{\tilde{Y}_1}}$, so $w(\alpha_1)$ belongs to the intersection $\overline{\mathcal{FE}_{\tilde{Y}_1}}\cap w(\overline{\mathcal{FE}_{\tilde{Y}_1}})$, and, by the same argument as in the proof of \cite[Theorem 1.3]{MarkmanLagrangian}, we obtain that $w$ is in the subgroup generated by reflections in divisors orthogonal to $\alpha_1$ inside $\mathcal{P}ex_{\tilde{Y}_1}$, i.e. $w(\alpha_1)=\alpha_1$. Finally, we get $h^*(\alpha_2)=\alpha_1$ and this shows the claim.
\eprf

It is worth remarking the following.

\begin{cor}
 Let $g\in\check{M}$ be an isotropic element. If there exist $w_1,w_2\in \Mo^2_{Hdg}(Y)$ such that $w_k(\check{\imath}(g))=c_1(\mathcal{L}_k)$, for $k=1,2$, and $\mathcal{L}_k$ induce rational Lagrangian fibrations on $Y$, then they induce birational rational Lagrangian fibrations.
 
 In particular,  on strictly $(\check{M},\check{\jmath})$-polarized manifolds, all the rational Lagrangian fibrations corresponding to $Z_g$ in the sense of Theorem \ref{mirrorLag} are birational.
\end{cor}
\prf It is an obvious application of Proposition \ref{Hdg-mon-Y} with $\gamma=w_1\circ w_2^{-1}$.
\eprf

\begin{rem}
In the case of smooth $K3$ surfaces, given $g$ an isotropic class in the boundary of the positive cone,  the sublattice $(\IZ g)^{\perp}/\IZ g$ has a clear geometrical meaning, as it corresponds to the configuration of the singular fibres. Here, we can only conjecture that, for non-birational rational Lagrangian fibrations induced by different $O(N)$-orbits of $g\in\check{M}$, the sublattices $(\IZ g)^{\perp}/\IZ g$ correspond to different configurations of divisors made of singular fibres. If the conjecture holds, then the results in this section would give not only an upper bound on the number of non-birational rational Lagrangian fibrations on a strictly $(\check{M},\check{\jmath})$-polarized IHS variety, but also indications on which singular fibres appear. This conjecture is the object of a forthcoming paper by the author and G. Sacc\`a.
\end{rem}

\section{Boundary components of $\mathcal{M}^+_{\langle 2d\rangle,j}$}\label{bound-comp}

In this section we are going to conclude the computation of the compactification in the case $M=\left\langle2d\right\rangle$ for fourfolds of $\hsk$-type: an application of the methods developed by Scattone in \cite{Scattone} will allow us to understand also one-dimensional boundary components.

Let $L=U^{\oplus 3}\oplus E_8^{\oplus 2}\oplus \left\langle-2\right\rangle$ be the $\hsk$ lattice.
Let $M\subset L$ be the rank one sublattice $\left\langle2d\right\rangle$ for a positive integer $d$. The following result, due to Gritsenko, Hulek and Sankaran, is known:

\begin{teo}\cite[Prop. 3.6 and 3.12]{GHS}
 The sublattice $M=\left\langle2d\right\rangle$ admits up to two non-isometric primitive embeddings in $L$. Let $h$ be a generator of $M$; then the following holds:
 \begin{enumerate}
  \item there is always a {\it split} embedding $j_s$, corresponding to $\divi h=1$, such that $N_s=U^{\oplus 2}\oplus E_8^{\oplus 2}\oplus \left\langle-2\right\rangle\oplus \left\langle-2d\right\rangle$, $\det N_s=4d$ and $A_{N_s}=\IZ/2\IZ\oplus \IZ/2d\IZ$;
  \item if $d\equiv 3$ modulo $4$, then $M$ admits a second embedding  $j_{ns}$, called {\it non-split}, corresponding to $\divi h=2$, such that $N_{ns}=U^{\oplus 2}\oplus E_8^{\oplus 2}\oplus B_d$, $\det N_{ns}=d$ and $A_{N_{ns}}=\IZ/d\IZ$.
 \end{enumerate}
In both cases, $\Gamma_{M,j}\cong \tilde{O}^+(N)$.
\end{teo}

Fix $j:M\subset L$ as above and let $\mathcal{M}_{2d}^+$ be the subset of $\mathcal{M}^+_{M,j}$, where $h$ corresponds to an ample class; as already shown in \cite{GHS}, we get an open algebraic embedding of $\mathcal{M}_{2d}^+/\Mo^2(M,j)$ into $D^+_M/\tilde{O}^+(N) $.

\begin{pro}
 The period map $\mathcal{P}_{M,j}^+$ restricts to an isomorphism $$\mathcal{M}_{2d}^+/\Mo^2(M,j)\rightarrow \left(D^+_{M}\setminus\coprod_{\delta\in\Delta(N)}(H_{\delta}\cap D^+_{M})\right)/\tilde{O}^+(N) $$
\end{pro}

In \cite[Lemma 5.6 and 5.7]{C4} we have explicitly computed $\tilde{O}(N)$-orbits of $m$-admissible isotropic elements in $N$ both in the split and in the non-split case, in order to understand how many non-isomorphic mirror families could occur in each case: 
each of these orbits gives a splitting of $N=U(m)\oplus \check{M}$ and hence a mirror family. Now we want to compute the number of boundary points in the Baily--Borel compactification: by \cite[Proposition 4.1.3]{Scattone}, we know that the zero-dimensional boundary components of $\overline{D^+_M/\tilde{O}^+(N)} $ are in bijection with  the set $I_1(A_N)/\lbrace \pm 1\rbrace$ of isotropic elements in $A_N$ modulo multiplication by $\pm 1$. 

We also want to compute the one-dimensional boundary components. Again, we apply the techniques developed by Scattone \cite{Scattone} in the $K3$ case, which rely on former work by Brieskorn \cite{Brieskorn}. The idea is the following: the set of one-dimensional rational boundary components is in bijection with the set $I_2(N)$ of rank two primitive isotropic sublattices of $N$ modulo the action of $\tilde{O}^+(N)$; to compute such orbits, we first compute $O(N)$-orbits and then we study the fibers of the surjective map $I_2(N)/\tilde{O}^+(N)\rightarrow I_2(N)/O(N)$.

\subsection{The split embedding}
Consider the primitive embedding $j_s:M=\left\langle2d\right\rangle\hookrightarrow L$, so that $N:=N_s=U^{\oplus 2}\oplus E_8^{\oplus 2}\oplus \left\langle-2\right\rangle\oplus \left\langle-2d\right\rangle$. In particular, if $e$ and $t$ denote respectively the generators of $\left\langle-2\right\rangle$ and $\left\langle-2d\right\rangle$ in $N_{s}$, the discriminant group $A_{N}$ is generated by $e/2$ and $t/2d$ and the discriminant quadratic form is given by 
\begin{displaymath}
q(\alpha \frac{t}{2d}+\beta \frac{e}{2})=-\frac{\alpha^2+\beta^2d}{2d}\in\IQ/2\IZ                                                                                                                                                                                                                                                                                                                                                                                                                                                                                 \end{displaymath}
for $\alpha=0,\dots,2d-1$ and $\beta=0,1$. 

Let $u$ and $v$ denote a standard basis of one of the two orthogonal summands $U$ inside $N_s$, and write $d=d'k^2$ with $d'$ square-free. We obtain the following:

\begin{pro}\label{zero-dim-split}
Zero-dimensional boundary components $F$ of $\overline{D^+_M/\tilde{O}^+(N)} $ are in bijection with $I_1(N)/\tilde{O}^+(N)$. Therefore their number is 
\[
\nu(d)=\left\lbrace
 \begin{array}{cc}
k+1&\mathrm{if}\ d'\equiv 3 (4)\\
&\\
\lfloor \frac{k+2}{2}\rfloor&\mathrm{otherwise}
\end{array}\right.
\]

\end{pro}
\prf
In the proof of \cite[Lemma 5.6]{C4} it was shown that isotropic elements in $A_N$ were either contained in $\IZ/2d\IZ$ or, if $d'\equiv 3 (4)$, of the form $\alpha(t/m+e/2)$ for $\alpha=0,\dots,m-1$, $m|2k$ but $m\nmid k$. We can rephrase that by saying that $I(A_N)=\cup_{m\mid K}H_m$, where:
\[
K=\left\lbrace
 \begin{array}{cc}
2k&\mathrm{if}\ d'\equiv 3 (4)\\
&\\
k&\mathrm{otherwise}
\end{array}\right.
,
\ H_m=\left\lbrace
 \begin{array}{cc}
\langle \frac{t}{m}\rangle&\mathrm{if}\ m|k\\
&\\
\langle \frac{t}{m}+\frac{e}{2}\rangle&\mathrm{if}\ m\nmid k,\ \frac{m}{2}|k\ \mathrm{and}\ d'\equiv 3 (4)
\end{array}\right.
\]
In particular, the isotropic elements of order $m$ are exactly
\[
x_{m,n}=\left\lbrace
 \begin{array}{cc}
\frac{nt}{m}&\mathrm{if}\ m|k,\\
\frac{nt}{m}+\frac{ne}{2}&\mathrm{if}\ m\nmid k,\ \frac{m}{2}|k\ \mathrm{and}\ d'\equiv 3 (4)
\end{array}\right.
\]
for $n=0,\dots,m-1$ with $\gcd(m,n)=1$.
The action of $\lbrace \pm 1\rbrace$ then interchanges $x_{m,n}$ and $x_{m,m-n}$, so that we can choose as representatives of the orbits the elements $x_{m,n}$ as above for $0\leq n\leq\frac{m}{2}$ with $\gcd(m,n)=1$.

Moreover, the number of zero-dimensional boundary components equals the order of $H_K/\lbrace \pm 1\rbrace$.
\eprf

We now want to compute $I_2(N)/\Gamma_{M,j}$; we start by computing $O(N)$-orbits of primitive isotropic rank two sublattices $S\subset N$, by applying Brieskorn's and Scattone's method. Given $S\in I_2(N)$ primitive, we define the finite quotient $H_S:=((S^{\perp})^{\perp}_{N^*})/S\subset A_N$ (see \cite[page 76]{Brieskorn}), where both the orthogonals are taken inside $N^*$. Since $S$ is isotropic, $H_S$ is isotropic as well inside $A_N$; Proposition \ref{zero-dim-split} implies the following

\begin{cor}\label{isot-HS}
 All isotropic subgroups of $A_N$ are cyclic, of the form $H_m$ for some $m|K$.
\end{cor}

We consider then the partition $I_2(N)=\cup I_{2,m}(N)$, where $I_{2,m}(N)$ is the subset consisting of $S\subset N$ with $H_S\cong\IZ/m\IZ$, and $m|k$ or, if $d'\equiv 3 (4)$, $m\nmid k$ and $\frac{m}{2}|k$; such a partition is $O(N)$-invariant.

For any $S\in I_2(N)$, let $E$ be the quotient $S^{\perp}/S$, that has signature $(0,18)$.

\begin{pro}\cite[\S 4.1, Lemma pag.77]{Brieskorn}\label{discriminantE}
 There is an isomorphism between $A_E$ and the quotient $H_S^{\perp}/H_S$ endowed with the quadratic form induced by the discriminant form of $A_N$.
\end{pro}

\begin{cor}
Let $S\in I_{2,m}(N)$ for $m|K$ and let $E$ be the quotient $S^{\perp}/S$; let $\overline{q}$ be the discriminant quadratic form of $E$.
\begin{enumerate}
 \item If $m|k$, the discriminant group $A_E$ is isomorphic to $\IZ/\frac{2d}{m^2}\IZ\oplus\IZ/2\IZ$, with quadratic form matrix $\overline{q}=\mathrm{diag}(-\frac{m^2}{2d},-\frac{1}{2})$.
 \item If $d'\equiv 3 (4)$, $m\nmid k$ and $\frac{m}{2}|k$, the discriminant group $A_E$ is isomorphic to $\IZ/\frac{4d}{m^2}\IZ$, with quadratic form matrix $\overline{q}=(-\frac{m^2+4d}{8d})$.
\end{enumerate}

\end{cor}
\prf
Take $S\in I_{2,m}(N)$ and suppose first that $m|k$; then, $H_S\cong H_m$ and it is easily seen that $H_S^{\perp}\cong\langle\frac{mt}{2d}\rangle\oplus\langle\frac{e}{2}\rangle\cong \IZ/\frac{2d}{m}\IZ\oplus\IZ/2\IZ$. As a consequence of Proposition \ref{discriminantE}, we get $A_E\cong \IZ/\frac{2d}{m^2}\IZ\oplus\IZ/2\IZ$, and the induced quadratic form can be represented by the matrix $\mathrm{diag}(-\frac{m^2}{2d},-\frac{1}{2})$.

If $d'\equiv 3 (4)$, $m\nmid k$ and $\frac{m}{2}|k$, then $H_S^{\perp}\cong\langle\frac{mt}{4d}+\frac{e}{2}\rangle\cong \IZ/\frac{4d}{m}\IZ$. In this case we obtain
$A_E\cong \IZ/\frac{4d}{m^2}\IZ$, with the induced quadratic form having value $-\frac{m^2+4d}{8d}$ modulo $2\IZ$ on the generator $\frac{mt}{4d}+\frac{e}{2}$. 
\eprf 

In both cases, $\det E=\frac{4d}{m^2}$. We now apply Brieskorn's theory of normal forms.

\begin{lem}\label{normalform}
 Let $S\in I_{2,m}(N)$ for $m|K$ and suppose that $\gcd(m,\det E)=1$; there exists $\lbrace v_1,\dots,v_{22}\rbrace$ a basis of the lattice $N$ such that $S=\langle v_1,v_2\rangle$, $S^{\perp}=\langle v_1,\dots,v_{20}\rangle$ and the quadratic form in this basis is of the form
 \[Q=\left(
 \begin{array}{ccc}
  0&0&A\\
  0&B&C\\
  ^tA&^tC&D
 \end{array}
\right)
 \]
 where $$A=\left(
 \begin{array}{cc}
  0&1\\
  m&0  
 \end{array}
\right),\ D=\left(
 \begin{array}{cc}
  2\delta&0\\
  0&0  
 \end{array}
\right),\ \mathrm{with}\ 0\leq\delta< m,$$ 
$C=0$ and $B$ represents the quadratic form on $E$.
\end{lem}
\prf We can apply \cite[Lemma 5.4.1]{Brieskorn} and we obtain a basis $\lbrace v_1,\dots,v_{22}\rbrace$  of $N$ which satisfies everything stated above except for the specific forms of the matrices $A,C,D$.

Denote by $M_{18,2}(\IZ)$ the set of $18\times2$ matrices with integer coefficients.
By definition, $A$ is the matrix
$A=\left(
 \begin{array}{cc}
  0&a_1\\
  a_2&0  
 \end{array}
\right)$,
where $H_S\cong\IZ/a_1\IZ\oplus\IZ/a_2\IZ$ and $a_1|a_2$, $C$ is a representative of a class in the quotient $$R:= M_{18,2}(\IZ)/(B\cdot M_{18,2}(\IZ)+M_{18,2}(\IZ)\cdot A),$$ and $D=\left(
 \begin{array}{cc}
  2d_{11}&d_{12}\\
  d_{12}&2d_{22} 
 \end{array}
\right)$ with $0\leq d_{ij}<a_{3-j}$ for $j=1,2$.

In our case, the same proof of \cite[Lemma 5.2.1]{Scattone} works and shows that we can assume $a_1=1$ and $a_2=m$, $C=0$ under the assumption $\gcd(m,\det E)=1$, and $D=\left(
 \begin{array}{cc}
  2\delta&0\\
  0&0 
 \end{array}
\right)$ with $0\leq \delta<m$.
\eprf

Denote $T(m,\delta)$ the lattice of rank two with quadratic form $\left(
 \begin{array}{cc}
  0&m\\
  m&2\delta 
 \end{array}
\right)$.

\begin{cor}\label{splitting}
Let $S\in I_{2,m}(N)$ for $m|K$ and suppose that $\gcd(m,\frac{4d}{m^2})=1$; there is a splitting $N\cong U\oplus T(m,\delta)\oplus E$.
\end{cor}

From this information one deduces bounds on the number of elements in $I_2(N)/O(N)$. Let $\mathcal{T}_d(m)$ be the set of integers $0\leq\delta<m$ for which the quadratic form $q_{m\delta}$ on $A_{T(m,\delta)}$ is compatible with the splitting $q_{A_N}=q_{m\delta}\oplus \overline{q}$, which comes from Corollary \ref{splitting}. Thus far we have proved the following:

\begin{pro}\label{rk-two-o(n)}
If $m|K$ and $\gcd(m,\frac{4d}{m^2})=1$, the map $(T(m,\delta),E)\mapsto S$ gives a surjection 
$$\mathcal{T}_d(m)\times\mathcal{G}(0,18,\overline{q})\rightarrow I_2(N)/O(N).$$
\end{pro}
\prf Indeed, given a pair $(T(m,\delta),E)$, Corollary \ref{splitting} tells that there exists an isomorphism $\psi: U\oplus T(m,\delta)\oplus E\rightarrow N$. If  $w_1,w_2$ and $w_3,w_4$ are respectively bases of $T(m,\delta)$ and of $U$ for which the quadratic forms are the usual ones, then $S:=\psi(\IZ w_1+\IZ w_3)$ is an isotropic sublattice of $N$ of rank two, determined up to the action of $O(N)$. Surjectivity is clear from above.
\eprf

Studying injectivity would require a longer analysis. From now on, we restrict ourselves to the case of $d$ a square-free integer, so that $k=1$; it will be more than enough for the scope of our application.
\begin{cor}\label{rk-two-o(n)-k=1}
 If $d$ is square-free, there is a bijection $$\mathcal{G}(E_8^{\oplus 2}\oplus\langle -2\rangle\oplus\langle -2d\rangle)\rightarrow I_2(N)/O(N).$$
\end{cor}
\prf When $k=1$ the only possibility is $m=1$ and henceforth $\delta=0$. The only element in $\mathcal{T}_d(1)$ is the hyperbolic lattice $U$, the genus $\mathcal{G}(0,18,\overline{q})$ is exactly the genus of $E_8^{\oplus 2}\oplus\langle -2\rangle\oplus\langle -2d\rangle$  and each element in it uniquely determines an $O(N)$-orbit of $S\in I_2(N)$.
\eprf

For the next step, Scattone's work holds true with small modifications: we recall that $O^+(N)$ has index two inside $O(N)$; let $\pi_1:I_2(N)/O^+(N)\ra I_2(N)/O(N)$ be the projection and $S\in I_2(N)$ be an isotropic sublattice. We have the following exact sequence involving the stabilizers:
\[\xymatrix{
 1\ar@{->}[r]&\stab_{O^+(N)}(S)\ar@{->}[r]&\stab_{O(N)}(S)\ar@{->}[r]&O(N)/O^+(N)\ar@{->}[r]&\coker\pi_1\ar@{->}[r]&1}\]
The fibre of $\pi_1$ over a point $S$ is isomorphic to the quotient of the two orbits $O^+(N)\cdot S/O(N)\cdot S$, and hence also to $\coker\pi_1$.

It follows from Lemma \ref{normalform} that
\begin{equation}\label{stabfiber}
\stab_{O(N)}(S)=\left\lbrace g\in GL(22,\IZ)|^tgQg=Q,\ 
g=\left(
\begin{array}{ccc}
U&V&UW\\
0&X&Y\\
0&0&Z
\end{array}\right)\right\rbrace
\end{equation}

It will also be useful, in order to study and compare the stabilizers with respect to different subgroups, to consider the map 
\[r:g\in GL(22,\IZ)\mapsto U\in GL(2,\IZ);\]
let $N_{\Gamma}(S)$ be $r(\stab_{\Gamma}(S))\subset GL(2,\IZ)$ for $\Gamma\subset O(N)$ an arithmetic subgroup.

\begin{lem}\label{NGammaS}
For any $S\in I_{2}(N)$, we have:
\begin{enumerate}
\item $N_{O(N)}(S)= GL(2,\IZ)$;
\item $N_{O^+(N)}(S)=SL(2,\IZ)$;
\item $N_{\tilde{O}^+(N)}(S)=SL(2,\IZ)$.
\end{enumerate}
\end{lem}
\prf The proof of \cite[Proposition 5.5.3 and Lemma 5.6.3]{Scattone} carries through, since it does not depend neither on the specific rank of $E$ nor on $A_N$.
\eprf

This is enough to conclude that $\stab_{O^+(N)}(S)$ and $\stab_{O(N)}(S)$ are not isomorphic, hence from (\ref{stabfiber}) we deduce that $\coker \pi_1$ is trivial.
\begin{cor}\label{O(N)=O+(N)}
There is an isomorphism $I_2(N)/O(N)\cong I_2(N)/O^+(N)$.
\end{cor}

Our problem is now reduced to the study of the fibres of the map $$I_2(N)/O^+(N)\ra I_2(N)/\tilde{O}^+(N)$$ We need to compute the cokernel of $$\pi_2: \stab_{O^+(N)}(S)\ra O^+(N)/\tilde{O}^+(N)\cong O(A_N),$$
hence we look at the quotient $\stab_{O^+(N)}(S)/\stab_{\tilde{O}^+(N)}(S)$.

Given $g\in\stab_{O^+(N)}(S)$ as in (\ref{stabfiber}), let $\alpha(g)$ be the class of $r(g)$ modulo $N_{\tilde{O}^+(N)}(S)$, and let $\xi(g)=\tau(X)\in O(A_E)$ for $$\tau:O(E)\ra O(A_E)$$ the standard projection map. Obviously, $\ker\alpha\supset \stab_{\tilde{O}^+(N)}(S)$; we denote $\overline{\alpha}$ the induced morphism $\stab_{O^+(N)}(S)/\stab_{\tilde{O}^+(N)}(S)\ra N_{O^+(N)}(S)/N_{\tilde{O}^+(N)}(S)$.
For a subgroup $\Gamma\subset O(N)$, denote $$\Fix_{\Gamma}(S)=\lbrace g\in O(N)|g_{|S}=\id_S\rbrace;$$ it fits into a short exact sequence of the following form:
\[
\xymatrix{
1\ar@{->}[r]&\Fix_{\Gamma}(S)\ar@{->}[r]&\stab_{\Gamma}(S)\ar@{->}[r]&N_{\Gamma}(S)\ar@{->}[r]&1}
\]
By choosing once $\Gamma=\tilde{O}^+(N)$ and once $\Gamma=O^+(N)$, we see that $\ker\overline{\alpha}$ is exactly the quotient $\Fix_{O^+(N)}(S)/\Fix_{\tilde{O}^+(N)}(S)$. By Corollary \ref{splitting} applied with $m=1$, $\delta=0$ and $T(1,0)\cong U$, the projection $\beta:g\in\Fix_{O^+(N)}(S)\mapsto X\in O(E)$ restricts to $\tilde{\beta}:g\in\Fix_{\tilde{O}^+(N)}(S)\mapsto X\in \tilde{O}(E)$; moreover, if $g\in\ker\beta$, then $X=I_E$ and $g_{|A_E}=\id_{A_E}$,  hence also $g_{|A_N}=\id_{A_N}$ and  $\ker\beta\cong\ker\tilde{\beta}$. Finally we obtain $\Fix_{O^+(N)}(S)/\Fix_{\tilde{O}^+(N)}(S)\cong O(E)/\tilde{O}(E)$. 

By Lemma \ref{NGammaS} we get $N_{O^+(N)}(S)\cong N_{\tilde{O}^+(N)}(S)$, hence we can construct the following commutative diagram:
 \[
\xymatrix{
1\ar@{->}[r]&O(E)/\tilde{O}(E)\ar@{->}[r]^{\tau}\ar@{->}[d]^{\cong}&O(A_E)\ar@{->}[r]\ar@{->}[d]&O(A_E)/\mathrm{Im}\tau\ar@{->}[r]\ar@{->}[d]^{\cong}&1\\
1\ar@{->}[r]&\stab_{O^+(N)}(S)/\stab_{\tilde{O}^+(N)}(S)\ar@{->}[r]&O(A_N)\ar@{->}[r]&\coker \pi_2\ar@{->}[r]&1\\
 }
 \]

We have thus proven the following
\begin{teo}\label{one-dim-split}
If $d$ is a square-free positive integer, there is a bijection between the set $I_2(N)/\tilde{O}^+(N)$ and the set 
\[
\lbrace (E,l)|E\in\mathcal{G}(E_8^{\oplus 2}\oplus\langle -2\rangle\oplus\langle -2d\rangle),\ l\in O(A_E)/\mathrm{Im}\tau\rbrace.
\]
 
\end{teo}

\subsection{The non-split embedding}

Consider the primitive embedding $j_{ns}:M=\left\langle2d\right\rangle\hookrightarrow L$, so that $N:=N_{ns}=U^{\oplus 2}\oplus E_8^{\oplus 2}\oplus B_d $.  In this case $d\equiv 3\ (4)$. Let $e$ and $w_1$, $w_2$ denote respectively the generators of $\left\langle-2\right\rangle$ and of a copy of $U$ in $L$, so that $M=\left\langle h \right\rangle\subset U\oplus\left\langle -2\right\rangle$ with $h=2w_1+\frac{d+1}{2}w_2+e$. The orthogonal $B_d$ is generated by $b_1=w_1-\frac{d+1}{4}w_2$ and $b_2=w_2+e$. The discriminant group $A_{N_{ns}}$ is generated by $t=\frac{1}{d}h-w_2$, and the discriminant quadratic form $q_{ns}$ is given by 
\begin{displaymath}
q_{ns}(\alpha t)=-\frac{2\alpha^2}{d}\in\IQ/2\IZ                                                                                                                                                                                                                                                                                                                                                                                                                                                                                 \end{displaymath}
for $\alpha=0,\dots,d-1$. 

Let $u$ and $v$ denote a standard basis of one of the two orthogonal summands $U$ inside $N_{ns}$, and write $d=d'k^2$ with $d'\equiv 3\ (4)$ square-free.
 We obtain the following:

\begin{pro}\label{zero-dim-non-split}
Zero-dimensional boundary components $F$ of $\overline{D^+_M/\tilde{O}^+(N)} $ are in bijection with $I_1(N)/\tilde{O}^+(N)$. Therefore their number is 
\[
\nu(d)=\lfloor \frac{k+2}{2}\rfloor
\]

\end{pro}
\prf
In the proof of \cite[Lemma 5.6]{C4} it was shown that isotropic elements in $A_N$ were all of order $m|k$. We can rephrase that by saying that $I(A_N)=\cup_{m\mid k}H_m$, where
$
H_m=
\langle \frac{dt}{m}\rangle
$.
In particular, the isotropic elements of order $m$ are exactly
$
x_{m,n}=
\frac{ndt}{m}$
for $n=0,\dots,m-1$ with $\gcd(m,n)=1$.
The action of $\lbrace \pm 1\rbrace$ then interchanges $x_{m,n}$ and $x_{m,m-n}$, so that we can choose as representatives of the orbits the elements $x_{m,n}$ as above for $0\leq n\leq\frac{m}{2}$ with $\gcd(m,n)=1$.

Moreover, the number of zero-dimensional boundary components equals the order of $H_k/\lbrace \pm 1\rbrace$.
\eprf

We now want to compute $I_2(N)/\Gamma_{M,j}$; once more, we start by computing $O(N)$-orbits of primitive isotropic rank two sublattices $S\subset N$. In this case, the fact that all isotropic subgroups of $A_N$ are cyclic is trivial.

We consider the partition $I_2(N)=\cup I_{2,m}(N)$, where $I_{2,m}(N)$ is the subset consisting of $S\in I_2(N)$ with $H_S\cong\IZ/m\IZ$ and $m|k$.

\begin{cor}
Let $S\in I_{2,m}(N)$ for $m|k$ and let $E$ be the quotient $S^{\perp}/S$; let $\overline{q}$ be the discriminant quadratic form of $E$.

The discriminant group $A_E$ is isomorphic to $\IZ/\frac{d}{m^2}\IZ$, with quadratic form matrix $\overline{q}=(-\frac{2m^2}{d})$.
\end{cor}
\prf
Take $S\in I_{2,m}(N)$; then $H_S\cong H_m$, and it is easily seen that $H_S^{\perp}\cong\langle mt\rangle\cong \IZ/\frac{d}{m}\IZ$. As a consequence of Proposition \ref{discriminantE}, we get $A_E\cong \IZ/\frac{d}{m^2}\IZ$, with the induced quadratic form having value $-\frac{2m^2}{d}$ modulo $2\IZ$ on the generator $mt$. 
\eprf 

In particular, $\det E=\frac{d}{m^2}$. Moreover, Lemma \ref{normalform} holds also in this case.

\begin{cor}\label{splitting-non-split}
Let $S\in I_{2,m}(N)$ for $m|k$ and suppose that $\gcd(m,\frac{d}{m^2})=1$; there is a splitting $N\cong U\oplus T(m,\delta)\oplus E$.
\end{cor}

\begin{pro}\label{rk-two-o(n)-non-split}
If $m|k$ and $\gcd(m,\frac{d}{m^2})=1$, the map $(T(m,\delta),E)\mapsto S$ gives a surjection 
$$\mathcal{T}_d(m)\times\mathcal{G}(0,18,\overline{q})\rightarrow I_2(N)/O(N).$$
\end{pro}

From now on, we again restrict ourselves to the case of $d$ a square-free integer, so that $k=1$.
\begin{cor}\label{rk-two-o(n)-k=1-non-split}
 If $d$ is square-free, there is a bijection $$\mathcal{G}(E_8^{\oplus 2}\oplus B_d)\rightarrow I_2(N)/O(N).$$
\end{cor}

Lemma \ref{NGammaS} and Corollary \ref{O(N)=O+(N)} are true in the non-split case as well.

Our problem is again reduced to the study of the fibres of the map $$I_2(N)/O^+(N)\ra I_2(N)/\tilde{O}^+(N).$$ Exactly the same proof as above shows  the following
\begin{teo}\label{one-dim-non-split}
If $d$ is a square-free positive integer, there is a bijection between the set $I_2(N)/\tilde{O}^+(N)$ and the set 
\[
\lbrace (E,l)|E\in\mathcal{G}(E_8^{\oplus 2}\oplus B_d),\ l\in O(A_E)/\mathrm{Im}\tau\rbrace.
\]
 
\end{teo}

\subsection{An application}

In this section we want to compute explicitly the boundary components of the Baily--Borel compactification of the moduli space $\mathcal{M}^+_{\langle 2\rangle,j_s}$ of $(\langle 2 \rangle, j_s)$-polarized fourfolds of $\hsk$-type: this is the moduli space of the so-called double EPW sextics (see \cite{O'Grady} for the definition).

\begin{cor}
The Baily--Borel compactification of $D^+_{\langle 2\rangle}/\tilde{O}^+(N)$ has exactly one $0$-dimensional boundary component.
\end{cor}

The first step in our computation is thus the computation of the lattices $E$ in the genus of $E_8^{\oplus 2}\oplus \left\langle-2\right\rangle^{\oplus 2}$. One standard strategy to compute these is to look for possible primitive embeddings of $E$ into a Niemeier unimodular lattice $\mathcal{N}$ of rank $24$, associated to a root system $R(\mathcal{N})$: if such an embedding exists, the orthogonal has to be isomorphic to $D_6$ by Corollary \ref{rk-two-o(n)-k=1}. In this way, one can find the root sublattice $R(E)$ and information about its embedding. This computation has already been performed by Nishiyama while studying Jacobian fibrations on some singular $K3$ surfaces.

\begin{pro}\cite[Theorem 3.1]{Nishiyama}
There are $13$ non-isometric lattices $E$ inside the genus of $E_8^{\oplus 2}\oplus \left\langle-2\right\rangle^{\oplus 2}$. Table \ref{genustable} lists all the possibilities for the root system $R(E)$ and for the root system $R(\mathcal{N})$ of the corresponding Niemeier lattice in which $E$ can be embedded.
\end{pro}

It follows from Theorem \ref{one-dim-split} that, in order to obtain the order of $I_2(N)/\tilde{O}^+(N)$, it would still be necessary to check whether the map $\tau:O(E)\ra O(q_E)$ is surjective or not for all possible $E\in\mathcal{G}(E_8^{\oplus 2}\oplus \left\langle-2\right\rangle^{\oplus 2})$.

\begin{rem}
 In \cite{dolgachevmirror}, the author remarks that the four elements in the genus of $E_8^{\oplus2}\oplus A_1$ correspond to the possible configurations of reducible fibres in an elliptic fibration $\epsilon$ on a strictly $(U\oplus E_8^{\oplus2}\oplus A_1)$-polarized $K3$ surface $\Sigma$. Four of the lattices $E$ in Table \ref{genustable} are exactly the direct sum of one element of the genus $E_8^{\oplus2}\oplus A_1$ and of an $A_1$ summand, namely $E_8^{\oplus 2}\oplus A_1 ^{\oplus 2}$, $D_{16}\oplus A_1 ^{\oplus 2}$, $E_7\oplus D_{10}\oplus A_1$ and $A_{17}\oplus A_1$.
 
 Given  $\Sigma$ an elliptic strictly $(U\oplus E_8^{\oplus2}\oplus A_1)$-polarized $K3$ surface, the IHS fourfold $\Sigma^{[2]}$ is strictly $(U\oplus E_8^{\oplus2}\oplus A_1^{\oplus2},\check{j})$-polarized and it naturally admits a Lagrangian fibration $f^{[2]}:\Sigma^{[2]}\rightarrow\IP^2$. Indeed, if $f:\Sigma\rightarrow \IP^1$ is an elliptic fibration, we get an induced fibration $f^{(2)}$ of the second symmetric product $\Sigma^{(2)} $ over $(\IP^1)^{(2)}\cong\IP^2$. The so-called {\it natural Lagrangian fibration} $f^{[2]}$ is the composition $f^{(2)}\circ \pi$, where $\pi$ is the resolution $\Sigma^{[2]}\ra \Sigma^{(2)}$. Let $\Delta_{\Sigma}\subset \Sigma^2$ be the diagonal and let $F$ be the exceptional divisor on $\Sigma^{[2]}$. By \cite[Theorem 1]{Fu}, $f^{[2]}(F)$ is one-dimensional, so that the generic fiber $(f^{[2]})^{-1}(p,q)$, for $(p,q)\in (\IP^1)^{[2]}\cong\IP^2$ not in $f^{[2]}(F)$, does not intersect $F$ and is isomorphic to the generic fiber of $f^{(2)}$, i.e. it is isomorphic to the product of the two elliptic curves $f^{-1}(p)$ and $f^{-1}(q)$.
 On the other hand, if we look at natural Lagrangian fibrations on $\Sigma^{[2]}$ induced by the four elliptic fibrations $\epsilon$ above, we can study their singular fibres: the discriminant $\Delta_{\epsilon^{[2]}}$ of the fibration $\epsilon^{[2]}$ in $\IP^2\cong\mathrm{Sym}^2(\IP^1)$ is a union of rational curves $D_i$, obtained as the product $\lbrace q_i\rbrace\times\IP^1$ for $q_i\in\Delta_{\epsilon}$, and of the diagonal $\Delta_{\IP^2}$. Divisors of singular fibres are, as a consequence, respectively threefolds fibered over $D_i$ with fibres isomorphic to a chain of bielliptic surfaces (compare with \cite{Matsushita} and \cite{Oguiso-Hwang}), and the exceptional divisor over the diagonal.
 
 Conjecturally, these divisors give a geometrical realization of the sublattice $E$. In this example the number of Lagrangian fibrations on strictly $(\check{M},\check{\jmath})$-polarized fourfolds of $\hsk$-type, up to birational transformation, is included between four and the order of $I_2(N)/\tilde{O}^+(N)$.
\end{rem}

\begin{table}[ht]
 
\caption{Genus of $E_8^{\oplus 2}\oplus \left\langle-2\right\rangle^{\oplus 2}$}\label{genustable}
\begin{center}

\begin{tabular}{|l|l|l|l|l|l|}
 \hline
$R(E)$ & $E_8^{\oplus 2}\oplus A_1 ^{\oplus 2}$& $D_{16}\oplus A_1 ^{\oplus 2}$ & $E_8\oplus D_{10}$ & $E_7\oplus D_{10}\oplus A_1$ & $E_7^{\oplus 2}\oplus D_4$  \\ \hline
$R(\mathcal{N})$ & $E_8^{\oplus 3}$ & $E_8\oplus D_{16}$ & $E_8\oplus D_{16}$ & $E_7^{\oplus 2}\oplus D_{10}$ & $E_7^{\oplus 2}\oplus D_{10}$  \\ \hline
$R(E)$ &  $A_{17}\oplus A_1$ & $D_{18}$& $D_{12}\oplus D_6$ & $A_1^{\oplus 2}\oplus D_{8}^{\oplus 2}$ & $A_3\oplus A_{15}$   \\ \hline
$R(\mathcal{N})$ & $E_{7}\oplus A_{17}$& $D_{24}$ & $D_{12}^{\oplus 2}$ & $D_{8}^{\oplus 3}$ & $D_9\oplus A_{15}$   \\ \hline
$R(E)$ &$E_6\oplus A_{11}\oplus\langle -4\rangle$ & $D_6^{\oplus 3}$&$A_9^{\oplus 2}$& &  \\ \hline
$R(\mathcal{N})$ & $E_6\oplus D_7\oplus A_{11}$ & $D_6^{\oplus 4}$&$D_6\oplus A_9^{\oplus 2}$ && \\ \hline

\end{tabular}
\end{center}

\end{table}

\bibliographystyle{amsplain}
\bibliography{BiblioModuliLPHK}
\end{document}